\documentclass[11pt]{article}
\usepackage{tikz,amsmath}
\usetikzlibrary{calc}

\newtheorem{thm}{Theorem}[section]
\newtheorem{lem}[thm]{Lemma}

\newcommand{\p}{\mathcal{P}}
\newcommand{\e}{\varepsilon}

\def\pf{\noindent{\it Proof.} }
\def\qed{\nopagebreak\hfill{\rule{4pt}{7pt}}\medbreak}

\makeatletter \@addtoreset{equation}{section} \makeatother

\begin{document}

\begin{center}
{\Large\bf On Universal Tilers}
\end{center}

\begin{center}
David G. L. Wang\\[6pt]

Beijing International Center for Mathematical Research\\
Peking University, Beijing 100871, P. R. China\\
{\tt wgl@math.pku.edu.cn}
\end{center}

\begin{abstract}
A famous problem in discrete geometry
is to find all monohedral plane tilers,
which is still open to the best of our knowledge.
This paper concerns with one of its variants that
to determine all convex polyhedra whose every
cross-section tiles the plane.
We call such polyhedra universal tilers.
We obtain that a convex polyhedron is a universal tiler
only if it is a tetrahedron or a pentahedron.
\end{abstract}

\noindent\textbf{Keywords:}
cross-section, Euler's formula,
pentahedron, universal tiler

\noindent\textbf{2010 AMS Classification:} 05B45


\section{Introduction}

A monohedral tiler is a polygon that can cover the plane
by congruent repetitions without gaps or overlaps.
The problem of determining all monohedral tilers,
also called the problem of tessellation,
was brought anew into mathematical prominence
by Hilbert when he posed it as one of his
``Mathematische Probleme'', see Kershner~\cite{Ker68}.
It is well-known that all triangles and
all quadrangles are tilers.
Reinhardt~\cite{Rei18} determined all
hexagonal tilers, and obtained some special kinds of pentagonal tilers.
Later it is shown that any
polygon with at least $7$ edges is not a tiler
by using Euler's formula,
see Dress and Huson~\cite{DH87}.
The problem of plane tiling,
however,
is now still open to the best of our knowledge.
In fact,
there are $14$ classes of pentagonal tilers were found,
see Hirschhorn and Hunt~\cite{HH85},
Sugimoto and Ogawa~\cite{SO06}, and Wells~\cite{Wel91}.
For a whole theory of tessellation patterns,
see Gr\"unbaum and Shephard's book~\cite{GS87} as a survey
up to 1987.

Considering a variant of the problem of plane tiling,
Akiyama \cite{Aki07} found all convex polyhedra
whose every development tiles the plane.
He call them tile-makers.
The main idea in his proof is
to investigate the polyhedra whose
facets tile the plane by stamping.
Notice that facets are special cross-sections.
This motivates us to
consider a more general class of polyhedron tilers.

Let $\p$ be a convex polyhedron,
and $\pi$ a plane.
Denote by $C(\pi)$ the intersection of $\pi$ and $\p$.
We say that
$\pi$ intersects $\p$ trivially
if $C(\pi)$ is empty, or a point,
or a line segment.
Otherwise we say $\pi$ intersects $\p$ non-trivially.
In this case,
$C(\pi)$ is a polygon with at least $3$ edges.
We call $C(\pi)$ a cross-section
if $\pi$ crosses $\p$ nontrivially.
We say that $\p$ is a universal tiler
if every cross-section of $\p$ tiles the plane.
In this paper,
we study the shape of universal tilers.
It is a variant of the problem of plane tiling.

It is easy to see that
every tetrahedron is a universal tiler
since every cross-section of a tetrahedron
is either a triangle or a quadrangle.
The main goal of this paper is to
present that any universal tiler has at most
$5$ facets.

This paper is organized as follows.
In Section~\ref{sec2},
we derive a necessary condition that a
hexagonal cross-section (if exists)
of a universal tiler satisfies.
It will be used for excluding the membership of
many polyhedra from the class of universal tilers.
In Section~\ref{sec3},
we prove that any facet of a universal tiler is either
a triangle or a quadrangle.
In Section~\ref{sec4},
by using Euler's formula
we obtain that any universal tiler has
at most $5$ facets.

\section{Hexagonal cross-sections of universal tilers}
\label{sec2}

Note that no polygonal tiler has more than $6$ edges.
It follows that any cross-section of a universal tiler
has at most $6$ edges.
In particular,
any facet of a universal tiler has no more than $6$ edges.
In this section,
we shall obtain a necessary condition
for hexagonal cross-sections of a universal tiler.

Let~$\p$ be a polyhedron,
and let~$\pi$ be a plane which
crosses~$\p$ nontrivially.
Let~$l$~be a line belonging to~$\pi$
and let $\e>0$.
Denote by~$\pi_+$ (resp.~$\pi_-$) the plane obtained by
rotating~$\pi$ around~$l$ by the angle~$\e$ (resp.~$-\e$).
Set $\e\to0$.
It is clear that
either~$\pi_+$ or~$\pi_-$
crosses~$\p$ nontrivially.
Of course it is possible that both~$\pi_+$ and~$\pi_-$
crosses~$\p$ non-trivially. Write
\[
p(\pi;\,l;\,\e)=
\begin{cases}
\pi_+,&\textrm{if $\pi_+$ crosses $\p$ nontrivially};\\[5pt]
\pi_-,&\textrm{otherwise}.
\end{cases}
\]
Then $p(\pi;\,l;\,\e)$ is a plane crossing~$\p$ nontrivially.
Intuitively, for small~$\e$,
the plane $p(\pi;\,l;\,\e)$ is obtained
by rotating the plane~$\pi$ a little along~$l$.
For notational simplification, we rewrite
\[
C(\pi;\,l;\,\e)=C(\,p(\pi;\,l;\,\e)\,).
\]
By the continuity of a polyhedron,
we see that the cross-section $C(\pi;\,l;\,\e')$
is nontrivial for any $0<\e'<\e$.
Let $C$ be a cross-section of~$\p$.
We say that $C$ is proper if
none of its vertices is a vertex of~$\p$,
that is,
any vertex of a proper cross-section
lies in the interior of an edge of~$\p$.

\begin{lem}\label{lem_proper}
If $\p$ has a cross-section with $n$ vertices,
then $\p$ has a proper cross-section with at least $n$ vertices.
\end{lem}

\pf
Set up an $xyz$-coordinate system.
For any real number~$a$,
denote by $\pi^a$ the plane determined by the equation $z\!=\!a$.
Suppose that the cross-section $C(\pi^0)$ has~$n$ vertices.
Without loss of generality,
we can suppose that the half-space~$z>0$ has
non-empty intersection with~$\p$.
By the continuity of~$\p$,
there exists $\delta>0$ such that
for any $0<\e<\delta$,
the cross-section $C(\pi^\e)$
has at least~$n$ vertices.
Consider the $z$-coordinates of all vertices of~$\p$.
Let~$\eta$ be the minimum positive $z$-coordinate among.
Then the cross-section $C(\pi^{\eta/2})$
is a proper cross-section with at least~$n$ vertices.
This completes the proof.
\qed

As will be seen,
with aid of the above lemma
we may take improper cross-sections out of our consideration.
Let $C(\pi)=V_1V_2\cdots V_n$
be a proper cross-section of~$\p$.
It is clear that for any $1\le i\le n$,
there is a unique edge of~$\p$ which contains~$V_i$,
denoted $e_i$.

\begin{lem}\label{lem_angle}
Let $C(\pi)=V_1V_2\cdots V_n$ be a proper cross-section
with $V_i\in e_i$.
Then there exists $\delta>0$ such that for any $0<\e<\delta$,
\begin{itemize}
\item[(i)]
$C(\pi;\,V_1V_3;\,\e)$ is a proper cross-section
with exactly $n$ vertices
which belong to the edges
$e_1$, $e_2$, $\ldots$, $e_n$ respectively; \item[(ii)]
if $C(\pi;\,V_1V_3;\,\e)=U_1^\e U_2^\e\cdots U_n^\e$,
where $U_i^\e\in e_i$, $U_1^\e=V_1$, and $U_3^\e=V_3$, then
\[
\angle V_1U_2^\e V_3\ne\angle V_1V_2V_3.
\]
\end{itemize}
\end{lem}

\pf
Since $C(\pi)$ is proper, by continuity,
there exists $\delta_1>0$ such that Condition (i) holds
for any $0<\e<\delta_1$.
Suppose that
\[
C(\pi;\,V_1V_3;\,\e)=V_1\,U_2^\e \,V_3\,U_4^\e\,
U_5^\e\,\cdots\, U_n^\e,
\]
where $U_i^\e\in e_i$.
Let $T$ be the trace of the point~$U_2^\e$
as~$\e$ varies such that
\begin{equation}\label{eq1}
\angle V_1U_2^\e V_3=\angle V_1V_2V_3.
\end{equation}
Then $T$ is a sphere if $\angle V_1V_2V_3=\pi/2$,
while $T$ is an ellipsoid otherwise.
On the other hand,
the point $U_2^\e$ moves along $e_2$
by Condition (i).
So~$U_2^\e$~belongs to the intersection
of a sphere (or ellipsoid) and a line.
Such an intersection contains at most two points,
say $\e_1$ and $\e_2$.
Taking $\delta<\min\{\delta_1,\e_1,\e_2\}$,
we complete the proof.
\qed

We need Reinhardt's theorem~\cite{Rei18}
of the classification of hexagonal tilers.
Traditionally,
we use the concatenation of two points,
say, $AB$, to denote
both the line segment connecting~$A$ and~$B$,
and its length.

\begin{thm}[Reinhardt]\label{thm_Reinhardt}
Let $V_1V_2\cdots V_6$ be a hexagonal tiler.
Then one of the following three properties holds:
\begin{itemize}
\item[(i)]
$V_1+V_2+V_3=2\pi$ and $V_3V_4=V_6V_1$;
\item[(ii)]
$V_1+V_2+V_4=2\pi$, $V_2V_3=V_4V_5$ and $V_3V_4=V_6V_1$;
\item[(iii)]
$V_1=V_3=V_5=2\pi/3$, $V_2V_3=V_3V_4$, $V_4V_5=V_5V_6$ and $V_6V_1=V_1V_2$.
\end{itemize}
\end{thm}

Figure 1 illustrates the $3$ classes of hexagonal tilers.
See also Bollob\'{a}s~\cite{Bol63} and
Gardner~\cite{Gar75} for its proof.

\begin{tikzpicture}
\begin{scope}
\draw
(0,0)
--node[below]{$a$}(1.5,0)node[above left]{$A$}--
++(60:1.3)node[left=2pt]{$B$}--
++(140:1)node[below]{$C$}--node[above]{$d$}
++(180:1.5)
--++(250:.8)
--cycle;
\node[below=.5cm,xshift=1cm,text width=3.5cm,text badly ragged]
{$A+B+C=2\pi$, $a=d$.};
\end{scope}
\begin{scope}[xshift=5cm]
\draw
(0,0)
--node[below]{$a$}
(1.5,0)node[above left]{$A$}--
++(60:1.2)node[left=3pt]{$B$}--node[above right]{$c$}
++(130:1)
--node[above right]{$d$}
++(170:1.5)node[below]{$D$}--node[above left]{$e$}
++(220:1)
--cycle;
\node[below=.5cm,xshift=.85cm,text width=3.5cm,text badly ragged]
{$A+B+D=2\pi$, $a=d$, $c=e$.};
\end{scope}
\begin{scope}[xshift=10cm]
\draw
(0,0)coordinate(f)
--node[below]{$a$}
(1.5,0)coordinate(a)node[above left]{$A$}--node[right]{$b$}
++(60:1.5)coordinate(b)
--node[above]{$c$}
++(150:1.7)coordinate(c) node[below]{$C$}--node[above]{$d$}
++(210:1.7)coordinate(d)
;
\coordinate (x) at ($(0,0)!.5!(d)$);
\coordinate (e) at ($(x)!{.5/sin(60)}!90:(d)$);
\draw
(d)--node[left]{$e$}(e)node[above right]{$E$}
--node[below]{$f$}(f);
\node[below=.5cm,xshift=.8cm,text width=3.5cm,text badly ragged]
{$A=C=E={2\over3}\pi$, $a=b$, $c=d$, $e=f$.};
\end{scope}
\begin{scope}[yshift=-2cm]
\node[xshift=6cm,text centered]
{Figure 1. The $3$ classes of hexagonal tilers.};
\end{scope}
\end{tikzpicture}

Denote by $\mathcal{H}_\p$ the set of proper hexagonal cross-sections of $\p$.

\begin{thm}\label{thm_hexcs}
Any proper hexagonal cross-section of a universal tiler, if exists,
has a pair of opposite edges of the same length.
\end{thm}

\pf Let $\p$ be a universal tiler with $\mathcal{H}_\p\ne\emptyset$.
For any $C\in\mathcal{H}_\p$,
denote by $a(C)$ the number of angles of size $2\pi/3$ in $C$.
Let
\[
S=\{H\in\mathcal{H}_\p\,\colon\,
\mbox{any pair of opposite edges of $H$ has distinct lengths}\}.
\]
Suppose to the contrary that $S\ne\emptyset$.
Let $H=V_1V_2\cdots V_6\in S$ such that
\[
a(H)=\min\{a(C)\,\colon\,C\in S\}.
\]
By Theorem~\ref{thm_Reinhardt}, we have $a(H)\ge3$.

Without loss of generality, we can suppose that
\begin{equation}\label{eq2pi/3}
\angle V_1V_2V_3={2\pi\over3}.
\end{equation}
By Lemma~\ref{lem_angle},
there exists $\delta>0$
such that for any $0<\e<\delta$,
\[
C(H;\,V_1V_3;\,\e)
=U_1^\e\, U_2^\e\, U_3^\e\,
U_4^\e\, U_5^\e\, U_6^\e\in\mathcal{H}_\p,
\]
where $U_1^\e=V_1$,
$U_3^\e=V_3$,
$U_i^\e\in e_i$ and
\begin{equation}\label{ineq4}
\angle\,{V_1\,U_2^\e\, V_3}\ne {2\pi\over3}.
\end{equation}
On the other hand, by continuity,
there exists $0<\eta<\delta$ such that
for any $i$ mod $6$,
\begin{align}
\Bigl|\,U_i^\eta\, U_{i+1}^\eta-U_{i+3}^\eta\, U_{i+4}^\eta\,\Bigr|
&\ge{1\over2}\,\Bigl|\,V_i\,V_{i+1}-V_{i+3}\,V_{i+4}\,\Bigr|,
\label{ineq-edge}\\[5pt]
\Bigl|\,\angle U_i^\eta\, U_{i+1}^\eta\, U_{i+2}^\eta
-{2\pi\over3}\,\Bigr|
&\ge{1\over2}\,\Bigl|\,\angle V_i\,V_{i+1}\,V_{i+2}-{2\pi\over3}\,\Bigr|.\label{ineq-angle}
\end{align}
Write $H^\eta=C(H;\,V_1V_3;\,\eta)$.
Then $H^\eta\in S$ by~\eqref{ineq-edge}.
In view of~\eqref{eq2pi/3}, \eqref{ineq4} and~\eqref{ineq-angle},
we deduce that
$a(H^\eta)\le a(H)-1$,
contradicting to the choice of~$H$.
This completes the proof.
\qed

As will be seen,
we shall obtain that any universal tiler
has no hexagonal cross-sections.
But we need Theorem~\ref{thm_hexcs}
to derive this result.

\section{The valence-sets of universal tilers}\label{sec3}

In this section, we show that any facet of a universal tiler
is either a triangle or a quadrangle.
Let $F=V_1V_2\cdots V_n$ be a facet of $\p$.
Let $d_i$ be the valence of~$V_i$.
We say that the multiset $\{d_1,d_2,\ldots,d_n\}$
is the valence-set of $F$.
For example, the valence-set of any facet of a~tetrahedron
is $\{3,3,3\}$.

\begin{lem}\label{lem_cutQ}
Let $\p$ be a universal tiler.
Let $\{d_1,d_2,\ldots,d_n\}$ be
a valence-set of a facet of~$\p$.
Then for any $1\le h\le n$,
there is a cross-section of~$\p$ with
${\sum_{i=1}^nd_i}-d_h-2n+4$
edges. Consequently, we have
\begin{equation}\label{2n+2}
{\sum_{i=1}^nd_i}-d_h\le 2n+2.
\end{equation}
\end{lem}

\pf
Let $F=V_1V_2\cdots V_n$ be a facet of $\p$,
where~$V_i$ has valence~$d_i$.
It suffices to show for the case $h=1$.
We shall prove by construction.

For convenience,
we set up an $xyz$-coordinate system as follows.
First,
choose a point~$U_1$
from the interior of the edge $V_nV_1$.
Set~$U_1$ to be the origin.
Next,
choose~$U_2$ from the interior of $V_1V_2$,
and build the $x$-axis by putting $U_2$ on the positive $x$-axis.
Then,
build the $y$-axis such that $F$ lies on the $xy$-plane
and the $y$-coordinate of $V_1$ is negative.
Consequently,
all the other vertices $V_2,\ldots,V_n$
have positive $y$-coordinates.
Since $F$ is a facet, the convex polyhedron $\p$
must lie entirely in one of the two half-spaces
divided by the $xy$-plane.
Build the $z$-axis such that all points in~$\p$ have nonnegative $z$-coordinates.
Now we have an $xyz$-coordinate system.

Let $S=\{F'\ |\ \mbox{$F'$ is a facet of $\p$},\
F'\cap F\ne\emptyset,\
F'\ne F\}$ with $|S|=s$. It is easy to see that
\begin{equation}\label{s}
s=\sum_{i=1}^nd_i-2n.
\end{equation}
By continuity,
there exists $\delta>0$ such that
for any $0<\e<\delta$,
the cross-section $C(z\!=\!\e)$ has exactly $s$~vertices.
Here, as usual, the equation $z\!=\!\e$
represents the plane parallel to the $xy$-plane
with distance $\e$.
Write
\[
C(z\!=\!\e)=C_1^\e\, C_2^\e\,\cdots\, C_s^\e\,.
\]
Then for any vertex $C_j^\e$,
there is a unique vertex $V_i$ such that
$V_i$ and $C_j^\e$ lie in the same edge of~$\p$.
Denote this~$V_i$ by~$R_j^\e$.
Clearly $R_j^\e$ is independent of~$\e$.
So we can omit the superscript~$\e$
and simply write~$R_j$.
Without loss of generality, we can suppose that
\[
R_1=R_2=\cdots=R_{\,t}=V_1,\quad R_{\,t+1}=V_2,\quad R_s=V_n,
\]
where
\begin{equation}\label{t}
t=d_1-2.
\end{equation}
Let $t+1\le k\le s$,
and let $y_k^\e$ be the $y$-coordinate of $C_k^\e$.
Since $V_2,\ldots,V_n$ have positive $y$-coordinates,
there exists $0<z_0<\delta$ such that
$y_k^\e>0$ for any $0<\e\le z_0$.
For simplifying notation,
we rewrite
\[
C(z\!=\!z_0)=C_1C_2\cdots C_s.
\]
Let $y_k$ be the $y$-coordinate of $C_k$. Set
\begin{equation}\label{epsilon0}
\e_0={1\over2}\min\Bigl\{\,z_0,\,{ z_0\over y_{t+1}},\,{ z_0\over y_{t+2}},\,
\ldots,\,{ z_0\over y_s}\,\Bigr\}.
\end{equation}

We shall show that
the cross-section $C(\pi_0)$ has
${\sum_{i=2}^nd_i}-2n+4$ edges.
Consider the function $f$ defined by
\[
f(V)=\e_0y-z,
\]
where $V=(x,y,z)$ is a point.
Denote by $\pi_0$ the plane determined by the equation $f(V)=0$.
On one hand, we have $f(R_k)>0$
since the vertex~$R_k$
has positive $y$-coordinate and zero $z$-coordinate.
On the other hand, by~\eqref{epsilon0} we have
\[
f(C_k)=\e_0y_k-z_0\le{1\over2}\cdot{z_0\over y_k}\cdot y_k-z_0<0.
\]
Therefore, the points~$R_k$ and~$C_k$
lie on distinct sides of~$\pi_0$.
Consequently, the plane $\pi_0$ intersects
the line segment~$C_kR_k$.
Let $I_k$ be the intersecting point.
Recall that $U_1$ is the origin
and $U_2$ lies on the positive $x$-axis.
So these two points belong to the plane~$\pi_0$.
Hence
\[
C(\pi_0)=U_1U_2I_{t+1}I_{t+2}\cdots I_s.
\]
By~\eqref{s} and~\eqref{t}, the number of edges of $C(\pi_0)$ is
\[
s-t+2=\sum_{i=1}^nd_i-2n-(d_1-2)+2={\sum_{i=2}^nd_i}-2n+4.
\]
Since any cross-section of a universal tiler
has at most $6$ edges, the inequality~\eqref{2n+2}
follows immediately.
This completes the proof.
\qed

\begin{lem}\label{lem_33333}
The valence-set of any facet of a universal tiler is not $\{3,3,3,3,3\}$.
\end{lem}

\pf
Let $\p$ be a universal tiler.
Suppose to the contrary that $\p$ has a pentagonal facet
$F$ whose every vertex has valence $3$.
For convenience, write $F=U_1'U_2U_3U_4U_5$.
Pick a point $U_1$ from the interior of the line segment~$U_1'U_2$
such that
\begin{equation}\label{ineq2}
U_1U_2\ne U_4U_5.
\end{equation}
Pick a point $U_6$ from the interior of the line segment~$U_5U_1'$
such that
\begin{equation}\label{ineq3}
U_2U_3\ne U_5U_6\quad\mbox{and}\quad U_3U_4\ne U_6U_1.
\end{equation}
The existences of~$U_1$ and~$U_6$ are clear.
Since the valence of each vertex of~$F$ is~$3$,
there exists~$\delta$ such that
for any $0<\e<\delta$,
\[
C(F;\,U_6U_1;\,\e)
=U_1^\e\, U_2^\e\, U_3^\e\, U_4^\e\, U_5^\e\, U_6^\e\in\mathcal{H}_\p,
\]
where $U_1^\e=U_1$, $U_6^\e=U_6$,
and $U_i^\e$ and $U_i$
lie on the same edge of $\p$ for each $2\le i\le 5$.
On the other hand, by continuity,
there exists $0<\eta<\delta$ such that
for any $i$ mod $6$,
\begin{equation}\label{u}
\bigl|\,U_i^\eta\, U_{i+1}^\eta
-U_{i+3}^\eta\, U_{i+4}^\eta\,\bigr|
\ge{1\over 2}\,\bigl|\,U_i\,U_{i+1}-U_{i+3}\,U_{i+4}\,\bigr|.
\end{equation}
In light of~\eqref{ineq2}, \eqref{ineq3} and~\eqref{u},
we see that the cross-section $C(F;\,U_1U_6;\,\eta)$ has
no pair of opposite edges of the same length,
contradicting to Theorem~\ref{thm_hexcs}.
This completes the proof.
\qed

Using similar combinatorial arguments as
in the above proof,
we can determine the shape of a facet of a universal tiler.

\begin{thm}\label{thm_facet}
Let $\p$ be a universal tiler.
Then every facet of $\p$ is either a triangle or a quadrangle.
Moreover,
the valence-set of any triangular facet (if exists)
of~$\p$ is either $\{4,3,3\}$ or $\{3,3,3\}$,
while the valence-set of any quadrilateral facet (if exists)
of~$\p$ is $\{3,3,3,3\}$.
\end{thm}

\pf
Let $F_n=V_1V_2\cdots V_n$ be a facet of $\p$.
Let $S_n=\{d_1,d_2,\ldots,d_n\}$
be the valence-set of~$F_n$.
By Lemma~\ref{lem_cutQ},
we see that for any $1\le h\le n$,
\[
2n+2\ge{\sum_{i=1}^nd_i}-d_h\ge3(n-1).
\]
Namely $n\le5$. If $n=5$, then~\eqref{2n+2} reads
\[
\sum_{i=1}^5d_i-d_h\le 12.
\]
Since each $d_i\ge3$, we deduce that
all $d_i=3$, contradicting to Lemma~\ref{lem_33333}.
Hence $n\le 4$.

If $n=4$, then
the valence-set~$S_4$ is
either $\{3,3,3,3\}$ or $\{4,3,3,3\}$ by~\eqref{2n+2}.
Assume that $S_4=\{4,3,3,3\}$,
where~$V_1$ has valence~$4$.
Pick a point~$A$ from the interior of the line segment~$V_1V_2$
such that $V_1A\ne V_3V_4$,
and a point~$B$ from the interior of~$V_2V_3$
such that $AB\ne V_4V_1$.
Similar to the proof of Lemma~\ref{lem_33333},
we can deduce that there exists~$\eta$
such that the cross-section $C(F_4;\,AB;\,\eta)$ belongs to~$\mathcal{H}_\p$,
and it has no pair of opposite edges of the same length,
contradicting to Theorem~\ref{thm_hexcs}.
Hence $S_4=\{3,3,3,3\}$.

Consider the case $n=3$.
By Lemma~\ref{lem_cutQ},
the valence-set~$S_3$
has five possibilities:
\[
\{3,3,3\},\quad\{4,3,3\},\quad\{4,4,3\},\quad\{4,4,4\},\quad\{5,3,3\}.
\]
If $S_3=\{4,4,3\}$ or $S_3=\{4,4,4\}$,
we can suppose that both~$V_1$ and~$V_2$ have valence~$4$.
Pick a point~$A$ from the interior of~$V_1V_2$,
and~$B$ from~$V_2V_3$
such that $AB\ne V_3V_1$.
Again, there exists~$\eta$
such that $C(F_3;\,AB;\,\eta)$
has no pair of opposite edges of the same length,
contradicting to Theorem~\ref{thm_hexcs}.
If $S_3=\{5,3,3\}$,
we can suppose that~$V_1$ has valence~$5$.
Pick~$A$ from the interior of~$V_1V_2$
such that $V_1A\ne V_3V_1$, and~$B$ from~$V_2V_3$.
By similar arguments, we get a contradiction to Theorem~\ref{thm_hexcs}.
Hence the valence-set~$S_3$ is either
$\{3,3,3\}$ or $\{4,3,3\}$.
This completes the proof.
\qed

\section{The shapes of universal tilers}\label{sec4}

In this section,
we show that every universal tiler
at at most $5$ facets.

Let~$\p$~be a universal tiler.
Let~$f$ (resp.~$v$,~$e$) be the total number of
facets (resp. vertices, edges) of~$\p$.
Euler's formula reads
\begin{equation}\label{Euler}
f+v=e+2.
\end{equation}
It is well-known that
there are two distinct topological types of pentahedra.
One is the quadrilateral-based pyramids,
which has the parameters
\[
(v,e,f)=(5,8,5);
\]
the other is pentahedra composed of two triangular bases
and three quadrilateral sides, which has
\[
(v,e,f)=(6,9,5).
\]

Let~$f_i$ be the number of facets of~$i$~edges in~$\p$.
Let~$v_i$ be the number of vertices of valence~$i$ in~$\p$.
By Theorem~\ref{thm_facet}, we have
\begin{equation}
f=f_3+f_4\quad\mbox{and}\quad
v=v_3+v_4.
\end{equation}
Here is the main result of this paper.

\begin{thm}\label{thm_main}
A convex polyhedron is a universal tiler only if
it is a tetrahedron or a pentahedron.
\end{thm}

\pf Let $\p$ be a universal tiler.
By Theorem~\ref{thm_facet},
every facet of $\p$ has at most $4$ edges
and every vertex of $\p$ has valence at most $4$.

First, we deduce some relations by double-counting.
Counting
the pairs $(e',f')$ where $f'$ is a facet of $\p$ and $e'$ is an
edge of $f'$, we find that
\begin{equation}\label{f34}
3f_3+4f_4=2e.
\end{equation}
Counting the pairs $(v',e')$ where $e'$ is
an edge of $\p$ and $v'$ is a vertex of $e'$, we obtain
\begin{equation}\label{v34}
3v_3+4v_4=2e.
\end{equation}
Combining the relations from~\eqref{Euler} to~\eqref{v34},
we deduce that
\[
(3f_3+4f_4)+(3v_3+4v_4)=4e=4(v+f-2)=4(v_3+v_4+f_3+f_4-2),
\]
namely
\begin{equation}\label{eq_f+v=8}
f_3+v_3=8.
\end{equation}
On the other hand,
taking the difference of~\eqref{f34} and~\eqref{v34} yields
\begin{equation}\label{Euler_diff}
4(f_4-v_4)=3(v_3-f_3).
\end{equation}
Now we count the pairs $(v',T)$,
where $T$~is a triangular facet of~$\p$
and~$v'$~is a vertex of~$T$ having valence~$4$.
By Theorem~\ref{thm_facet},
every triangular facet has at most one vertex of valence~$4$,
and every facet containing a vertex of valence~$4$ must be a triangle.
Therefore
\begin{equation}\label{ineq_f3v4}
4v_4\le f_3.
\end{equation}
By~\eqref{eq_f+v=8}, \eqref{Euler_diff} and~\eqref{ineq_f3v4},
we deduce that $f_3\le4$.
Note that~$f_3$ is an even number by~\eqref{f34}.
So $f_3\in\{0,2,4\}$.

If $f_3=4$, then
$v_3=4$ by~\eqref{eq_f+v=8},
and $f_4=v_4\le1$ by~\eqref{Euler_diff}
and~\eqref{ineq_f3v4}.
In this case, $\p$~is a tetrahedron if $f_4=0$,
and~$\p$ is a quadrilateral-based pyramid if $f_4=1$.

If $f_3=2$, then $v_3=6$ by~\eqref{eq_f+v=8},
$v_4=0$ by~\eqref{ineq_f3v4},
and consequently $f_4=3$ by~\eqref{Euler_diff}.
In this case, $\p$ is a pentahedron composed of
two triangular bases and three quadrilateral sides.

If $f_3=0$, then $v_3=8$, $v_4=0$, and $f_4=6$.
Thus $\p$ is a cube.
We shall show that it is impossible.
Denote
\[
\p=ABCD\mbox{-}EFGH.
\]
For convenience,
we set up an $xyz$-coordinate system
such that the plane $z\!=\!0$ coincides
with the plane $ACH$,
and the vertex~$D$ has negative $z$-coordinate.
Let~$z_B$ (resp.~$z_E$, $z_F$, $z_G$)
be the $z$-coordinate of~$B$ (resp.~$E$, $F$, $G$).
Since $\p$ is convex, all these $z$-coordinates are positive.
Write
\[
\delta={1\over2}\min\{z_B,\,z_E,\,z_F,\,z_G\}.
\]
Then the line segment~$AB$ intersects the plane $z\!=\!\e$.
Let~$A_1^\e$~be the intersecting point.
Similarly,
let~$A_2^\e$
(resp.~$C_1^\e$, $C_2^\e$,
$H_1^\e$, $H_2^\e$)
be the intersection of the plane $z=\e$
and the line segment~$AE$
(resp.~$BC$, $CG$, $GH$, $HE$).
So
\[
C(z\!=\!\e)
=A_2^\e\, A_1^\e\, C_1^\e\, C_2^\e\, H_1^\e\, H_2^\e\in\mathcal{H}_\p.
\]
\begin{center}
\begin{tikzpicture}
\begin{scope}
\path
coordinate (O) at (0,0)
coordinate (x) at (5.5,0)
coordinate (y) at (80:3.7)
coordinate (z) at (240:2.5);


\draw coordinate (A) at (240:2)
  let \p1=($(A)$) in
  (A)node[below=.5mm]{$A$}coordinate(A)
  --(3,\y1) node[below=.5mm]{$B$}coordinate (B)
  --(5,0)node[below right]{$C$}coordinate(C);
\draw
  (110:1.5)node[left]{$E$}coordinate (E)
  --(25:3)node[right]{$F$}coordinate (F)
  --(35:4.5)node[above=.5mm]{$G$}coordinate (G)
  --(80:2.5)node[above=1mm]{$H$}coordinate (H)
  --cycle;

\draw[dotted]
(O)node[below right]{$D$}--(A)
(O)--(C)
(O)--(H);

\draw
(A)--(E)
(B)--(F)
(C)--(G);

\coordinate (A1) at (intersection of A--B and x--z);
\coordinate (A2) at (intersection of A--E and y--z);
\node[below of=A1,yshift=0.6cm]{$A_1^\e$};
\node[left of=A2,xshift=0.6cm]{$A_2^\e$};
\coordinate (C1) at (intersection of B--C and x--z);
\coordinate (C2) at (intersection of C--G and x--y);
\node[below of=C1,yshift=0.6cm]{$C_1^\e$};
\node[right of=C2,xshift=-0.5cm]{$C_2^\e$};
\coordinate (H1) at (intersection of G--H and x--y);
\coordinate (H2) at (intersection of H--E and y--z);
\node[above of=H1,yshift=-0.7cm]{$H_1^\e$};
\node[left of=H2,xshift=0.5cm]{$H_2^\e$};

\draw[thick]
(A1)--(A2)--(H2)--(H1)--(C2)--(C1)--cycle;

\foreach \point in {A1,A2,C1,C2,H1,H2}
\fill[black,opacity=1] (\point) circle (1.5pt);
\end{scope}
\begin{scope}[yshift=-3cm]
\node[xshift=2cm,text centered]
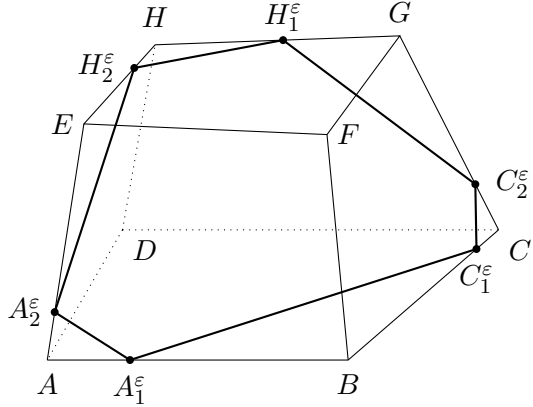
{Figure 2. The hexagonal cross-section $C(z\!=\!\e)
=A_2^\e\, A_1^\e\, C_1^\e\, C_2^\e\, H_1^\e\, H_2^\e$.};
\end{scope}
\end{tikzpicture}
\end{center}
By continuity, we have
\[
A_1^\e\, A_2^\e\,\to0,\quad
C_1^\e\, C_2^\e\,\to0,\quad
H_1^\e\, H_2^\e\,\to0,
\]
as $\e\to0$.
So there is $0<\eta<\delta$ such that the cross-section
$C(z\!=\!\eta)$ has no pair of opposite edges of the same length,
contradicting to Theorem~\ref{thm_hexcs}.
This completes the proof.
\qed

Recall that any tetrahedron $T$ is a universal tiler.
We present that pentahedron universal tilers also exist.

\begin{thm}\label{thm1}
Any pentahedron having a pair of parallel facets
is a universal tiler.
\end{thm}

\pf Suppose that $\p$ is a pentahedron with a pair of parallel facets.
Note that any cross-section of a pentahedron has at most $5$ edges.
It suffices to show that any pentagonal cross-section of $\p$
tiles the plane.
Let $C$ be a pentagonal cross-section of $\p$.
Then $C$ has a pair of parallel edges.
As pointed out by Reinhardt in~\cite{Rei18},
any pentagon with a pair of parallel edges is a tiler.
This completes the proof.
\qed

\end{document}